\documentclass[letterpaper, 10 pt, conference]{ieeetran}  

\IEEEoverridecommandlockouts                              

\usepackage{amsmath} 
\usepackage{amssymb}  

\usepackage{amsthm}

\newtheoremstyle{theoremdd}
{\topsep}
{\topsep}
{\itshape}
{0pt}
{\bfseries}
{:}
{ }
{\thmname{#1}\thmnumber{ #2}\thmnote{ (#3)}}

\theoremstyle{theoremdd}

\newtheorem{theorem}{Theorem}[section]

\newtheorem{lemma}[theorem]{Lemma}
\newtheorem*{main_suff_thm}{Sufficiency Theorem}
\newtheorem*{main_thm}{Necessity Theorem}

\renewenvironment{proof}{{\bf Proof:}}

\usepackage{amssymb}  
\usepackage{mathrsfs}
\usepackage{graphicx} 
\usepackage{epsfig} 
\usepackage{epstopdf}
\usepackage{color}
\usepackage{caption}
\usepackage{float}
\usepackage{mathtools}

\usepackage{lipsum}
\makeatletter
\def\footnoterule{\relax%
	\kern-5pt
	\hbox to \columnwidth{\hfill\vrule width .9\columnwidth height 0.4pt\hfill}
	\kern4.6pt}
\makeatother

\makeatletter
\newcommand{\figcaption}{\def\@captype{figure}\caption}
\newcommand{\tabcaption}{\def\@captype{table}\caption}

\usepackage{color,hyperref}
\definecolor{darkblue}{rgb}{0.0,0.0,0.3}
\hypersetup{colorlinks,breaklinks,linkcolor=darkblue,urlcolor=darkblue,anchorcolor=darkblue,citecolor=darkblue}

\makeatother
\title{\LARGE \bf  On Positive Solutions of a Delay Equation\\ Arising  When Trading in Financial Markets 
}

\author{\large Chung-Han Hsieh,$^{*}$  B. Ross Barmish,$^{**}$ and John A. Gubner$^{***}$
	\thanks{\hskip -10pt ${}^*$Chung-Han Hsieh is a graduate student working towards a Ph.D. degree in the Department of Electrical and Computer Engineering, University of Wisconsin, Madison, WI 53706. E-mail: \href{mailto: hsieh23@wisc.edu}{hsieh23@wisc.edu}.} 
	\thanks{\hskip -10pt ${}^{**}$B. Ross Barmish is a Research Professor in the 
		Department of Elec\-tric\-al and Computer  Engineering, Boston University, Boston, MA\ \ 02215.
		E-mail: \href{mailto: barmish@bu.edu}{barmish@bu.edu} 	}
	\thanks{\hskip -10pt ${}^{***}$John A. Gubner is a Professor in  the Department of Electrical and Computer Engineering, University of Wisconsin, Madison, WI 53706. \mbox{E-mail}: \href{mailto: john.gubner@wisc.edu.}{john.gubner@wisc.edu} 	  \vspace{3mm} }
}

\begin{document}

	\maketitle
	\thispagestyle{empty}
	\pagestyle{empty}
	
	\parindent = 0pt
	
	\begin{abstract}
		We consider a discrete-time, linear state equation with delay which arises as a model for a trader's account value when buying and selling a risky asset in a financial market. The state equation includes a nonnegative feedback gain~$\alpha$ and a sequence~$v(k)$ which models asset returns which are within known bounds but otherwise arbitrary. 
		We introduce two thresholds,~$\alpha_-$ and~$\alpha_+$, depending on these bounds, and prove that 
		for~\mbox{$\alpha < \alpha_-$}, state positivity is guaranteed for all time and all asset-return sequences; i.e., bankruptcy is ruled out and positive solutions of the state equation are continuable indefinitely.
		On the other hand, for~$\alpha > \alpha_+$, we show that there is always a sequence of asset returns for which the state fails to be positive for all time; i.e., along this sequence, bankruptcy occurs and the solution of the state equation ceases to be meaningful after some finite time. 
		Finally, this paper also includes a conjecture which says that for the ``gap" interval~\mbox{$\alpha_- \leq \alpha \leq \alpha_+,$}   state positivity is also guaranteed for all time. 
		Support for the conjecture, both theoretical and computational, is~provided.
	\end{abstract}

	
	%
	\vspace{6mm}
	\section{Introduction}
	\label{SECTION: INTRODUCTION}
	\vspace{-1mm}
	The motivation for this paper is derived from an emerging line of research involving the use of system-theoretic ideas to trade in financial markets; e.g., see~\mbox{\cite{Cover_Ordentlich_1996}--\cite{Hsieh_Gubner_Barmish_2018_CDC}}. Similar to previous work, in this paper, we operate in an idealized market with no transaction costs such as brokerage commission or fees and with perfect liquidity; i.e., there is no gap between the bid and ask prices, and the trader has the ability to buy or sell any number of shares, including fractions, at the market price.  These assumptions arise in the finance literature in the context of ``frictionless" markets; e.g., see~\cite{Merton_1992}. 
	
	\vspace{3mm}
	With the above providing the backdrop, this paper concentrates on a difference equation with delay and establishes conditions under which all solutions $X(k)$ are positive for all $k$. We refer to this as ``all-time positivity."
	Related to this work on all-time positivity are papers in the mathematics literature which deal with  difference equations with multiple delays and provide conditions under which solutions are either {eventually} positive or {eventually} negative; i.e., $X(k)$ has one sign for $k$ suitably large; e.g., see \cite{Erbe_Zhang_1989} and \cite{Berezansky_Braverman_2006} and their bibliographies. 
	As noted in Section~\ref{SECTION: Proofs}, conditions in the aforementioned literature under which eventual positivity and negativity fail can be viewed as a special case of our theorem which provides a necessary condition for all-time~positivity. 

	\vspace{5mm}
	\textbf{Problem Formulation:}
	To formulate the problem at hand, we use $v(k)$ to represent the unpredictable \textit{returns} of a risky asset such as a stock or a foreign currency at stage $k$.   
	Our state-equation model for the account value $X(k)$ includes a delay due to the fact that a trader's interactions with the market are not instantaneous. 
	Specifically, at stage~$k$,  we take $\alpha \geq 0$ to be a feedback gain representing the targeted percentage of a trader's account $X(k)$ to be invested in the risky asset. 
	Then, at stage~$k$, order transmission and execution delay are accounted for by the realized control $u(k)$ representing the dollar level of investment at $k$. We begin with $u(0) \doteq 0$, and for~\mbox{$k>0$}, 
	\[
	u(k) \doteq \alpha (1+v(k-1)) X(k-1)
	\]
	to account for the delay.
	Accordingly, the closed-loop state equation is
	\begin{align*}
	X(k+1) 	& = X(k) + u(k) v(k) \\[.5ex]
	& = X(k) + \alpha (1+v(k-1))X(k-1)\,v(k)
	\end{align*}
	with positive initial conditions 
	$$
	X(0)=X(1)= X_0>0.
	$$ 
	In the sequel, a time-varying sequence of risky asset returns
	\[
	v \doteq \{v(k)\}_{k=0}^\infty
	\]
	is called a {\it path}, and is said to be {\it admissible} if it stays within  known~bounds
	$$
	v_{\min} \leq v(k) \leq v_{\max}
	$$ 
	where 
	$$
	-1 < v_{\min} < 0< v_{\max} <\infty.
	$$
	The assumption~\mbox{$v_{\min} > -1$} excludes the case that the underlying asset price can reach zero. 
	We take $\cal V$ to be the set of all admissible paths and often emphasize the state dependence on~$v \in \mathcal{V}$ by writing~$X(v,k)$  instead of~$X(k)$.  Additionally, we take $\mathcal{V}^N$ to be the set of all~\mbox{$v=(v(0),v(1),\ldots,v(N-1))$} such that for~\mbox{$k=0,1, \ldots,N-1$},  $v(k)$ stays within the known bounds above. Elements of ${\cal V}^N$ are called \textit{admissible partial paths}, or simply admissible paths when there is no confusion. 
	
	\vspace{3mm}
	As is typical in control theory, it is convenient to eliminate the delay term in the state equation above and work with a two-state system. That is, defining the state vector 
	$$
	x(k) \doteq [X(k) \;\; X(k-1)]^T,
	$$
	we obtain the linear time-varying system
	\[
	x(k+1) = A(v,k)\, x(k)
	\]
	where  
	\[
	A(v,k) \doteq 
	\begin{bmatrix}
	1 & \alpha (1+v(k-1))v(k) \\ 
	1 & 0
	\end{bmatrix}.
	\]
	As mentioned previously, we work with the specific initial conditions~\mbox{$X(0)=X(1) = X_0>0$}. 
	Although our goal, state positivity, is the same as in existing positive system theory, for example, see~\cite{Farina_2002} and~\cite{Kaczorek_2006}, this body of work  
	is not in play because the matrix $A(v,k)$ can have a negative entry. 


	\vspace{5mm}
	\textbf{All-Time Positivity:}
	Although the solution to the state equation exists for all $k$,
	since bankruptcy precludes future trading, the analysis ceases to be meaningful once $X(v,k) \leq 0$. With this as motivation, the focal point in this paper is the issue of all-time positivity. 
	In a sense, we are addressing a question about existence and continuability of positive solutions for infinitely many stages.
	Indeed, for a given feedback gain~\mbox{$\alpha \geq 0$}, we say that the \textit{all-time positivity} condition holds~if 
	$$
	X(v,k)>0
	$$
	for all $v \in \mathcal{V}$ and all $k \geq 0.$ 

	\vspace{3mm}
	It is also worth mentioning that~$u(k)\geq 0$ is guaranteed when all-time positivity holds. In finance, the condition $u(k) \geq 0$ is interpreted to mean that the trader holds a long position and no \textit{short selling} occurs. 
	Finally, we mention that traders in financial markets often have their orders restricted by leverage constraints imposed by the broker. 
	That is, letting
	$$
	L(k) \doteq \frac{u(k)}{X(k)},
	$$
	a maximum allowed leverage $L_{\max}>0$ is specified, and the trader's account is securitized by a requirement that~\mbox{$L(k) \leq L_{\max}$}. 
	For markets involving stock, $L_{\max}\leq 2$ is rather typical, and for foreign currency trading, \mbox{$L_{\max}\leq 100$} can easily be the case. 
	Leverage imposes a restriction on~$\alpha$. 
	However, since our criteria and conjecture on all-time positivity apply for all $\alpha \geq 0$, leverage bounds are ignored since they have no effect on the analysis to follow. 
	
	
	\vspace{5mm}
	{\bf Plan for the Remainder of the Paper:} In Section~\ref{SECTION: Main Results}, we present our main results. To this end, the section is centered around two critical thresholds, $\alpha_-$ and  $\alpha_+$ with~\mbox{$\alpha_- < \alpha_+$}, which we define.
	We first provide a result, called the Sufficiency Theorem, which tells us that~\mbox{$\alpha < \alpha_-$} is sufficient for all-time positivity. 
	Our next result, called the Necessity Theorem, gives a necessary condition for all-time positivity. 
	Specifically, for~\mbox{$\alpha > \alpha_+$}, we prove that there is a sequence of asset returns, called the \textit{distinguished} path and denoted by~$v^*$, for which the state fails to be positive for all~$k$. 
	In Section~\ref{SECTION: Technical Results}, we state two preliminary technical results regarding the state~$X(v^*,k)$ along this path. Next,
	in Section~\ref{SECTION: Proofs}, the proofs of the preliminary and main results   are provided.  
	In Section~\ref{SECTION: Conjs and Supports Of the Conjectures}, we provide a conjecture which says that all-time positivity  is guaranteed for the ``gap" interval~\mbox{$\alpha_- \leq \alpha \leq \alpha_+ $}. 
	%
	The section includes both theoretical and computational support for the conjecture. 
	Finally,  in Section~\ref{SECTION: CONCLUSION AND FUTURE WORK}, some concluding remarks are given, and possible directions for future research are~indicated.

	\vspace{6mm}
	\section{Main Results}
	\label{SECTION: Main Results}
	\vspace{0mm}
	The main results to follow involve two critical thresholds,~$\alpha_-$ and~$\alpha_+$. 
	The first of these, $\alpha_{-}$, is motivated by considering~\mbox{$k=2$} and noting that
	\begin{align*}
	X(2) 		& = X(1) + \alpha (1+v(0))v(1) \, X(0)\\[.5ex]
	& \geq [1+\alpha(1+v_{\max})v_{\min}]X_0.
	\end{align*}
	This lower bound is positive if and only if
	$$
	\alpha <  \frac{1} {|v_{\min}|(1+v_{\max})} .
	$$
	To show $X(k)>0$ for all $k$ rather than just $k=2$, the theorem below, proved in Section~\ref{SECTION: Proofs}, requires the stronger assumption~that $\alpha < \alpha_-$, where
	\[
	\alpha_- \doteq \frac{1}{1+v_{\max}}.
	\]

	\vspace{5mm}
	\begin{main_suff_thm}\label{thm: ASP holds (SP Lemma)}
		The condition~\mbox{$
			0\leq \alpha < \alpha_-
			$} is sufficient for all-time positivity. That is, if $\alpha < \alpha_{-}$, given any admissible path $v\in \mathcal{V}$, it follows that $X(v,k)>0$ for all~$k$.
	\end{main_suff_thm}
	

	\vspace{5mm}
	\textbf{Necessary Condition for All-Time Positivity:} As mentioned in the introduction, our necessary condition for all-time positivity is motivated by studying the state equation in response
	to a {\it distinguished path\/} of returns~$v^*$. This path is
	defined by~$v^*(0)=v_{\max}$ and~\mbox{$v^*(k)=v_{\min}$} for~\mbox{$k \geq 1 $}.
	Along this path, since the first trade is executed at stage~$k=1$, the return~$v(0) = v_{\max}$ can be viewed as ``baiting" the trader with a large positive return and then, a worst-case scenario of sorts occurs because the account loses value on every subsequent trade. 
	To motivate the definition of the threshold $\alpha_+$ entering into our analysis of necessity, 
	let
	$$
	\alpha_s \doteq \frac{1}{4|v_{\min}|(1+v_{\min})}.
	$$
	Then with $\alpha > \alpha_s,$ consistent with the fact that the matrix 
	$$
	A(v^*,k) = \begin{bmatrix}
	1 & \alpha (1+v_{\min})v_{\min} \\ 
	1 & 0
	\end{bmatrix}
	$$
	has a pair of complex eigenvalues, the state $X(v^*,k)$ is oscillatory about zero. Hence the value of state can be negative.  This becomes a special case of the theorem to follow.  
	For~\mbox{$\alpha \leq \alpha_s$},
	the solution is nonoscillatory, and 
	our analysis shows that the state is negative for large $k$ when $v_{\max}>1+2v_{\min}$ and~$\alpha > \alpha^*$, where
	\[
	\alpha^* \doteq \frac{v_{\max}-v_{\min}}{|v_{\min}|(1+v_{\max})^2}.
	\] 
	The theorem below, whose proof is relegated to Section~\ref{SECTION: Proofs}, brings these ideas to~fruition. The threshold $\alpha_{+}$ defined next is readily verified to exceed $\alpha_{-}.$


	\vspace{5mm}
	\begin{main_thm}\label{thm: ASP failure}
		With 
		\[ 
		\alpha_+ \doteq \begin{cases} \displaystyle
		\alpha^* & \text{ if } v_{\max} > 1+2v_{\min}; \\[1ex]
		\displaystyle \alpha_s &\text{ if } v_{\max} \leq 1+2v_{\min},
		\end{cases}
		\]
		the condition $\alpha \leq \alpha_{+}$ is necessary for all-time positivity. Equivalently, if $\alpha > \alpha_{+},$ then there exists an admissible path~\mbox{$v \in \mathcal{V}$} such that
		$
		X(v,k) \leq 0
		$
		for some $k$.
	\end{main_thm}

	\vspace{5mm}
	\textbf{Graphical Depiction of Bounds:} 
	In Figure~\ref{fig: conjectured alpha_max}, the dependencies of~\mbox{$\alpha_-$} and $\alpha_+$ on $v_{\min}$ and $v_{\max}$ are displayed over the range
	\mbox{$
		-1<v_{\min} < 0 < v_{\max}=2.
		$}
	The lower surface (black) is obtained by using the formula for~$\alpha_-$, the red part of the upper surface is obtained by using the formula for~\mbox{$\alpha_{+}=\alpha_s$}, and the larger green part of the upper surface is obtained using the formula for~$\alpha_{+}=\alpha^*$. 
	
	\vspace{3mm}
	\begin{center}
		\graphicspath{{Figs/}}
		\includegraphics[scale=0.4]{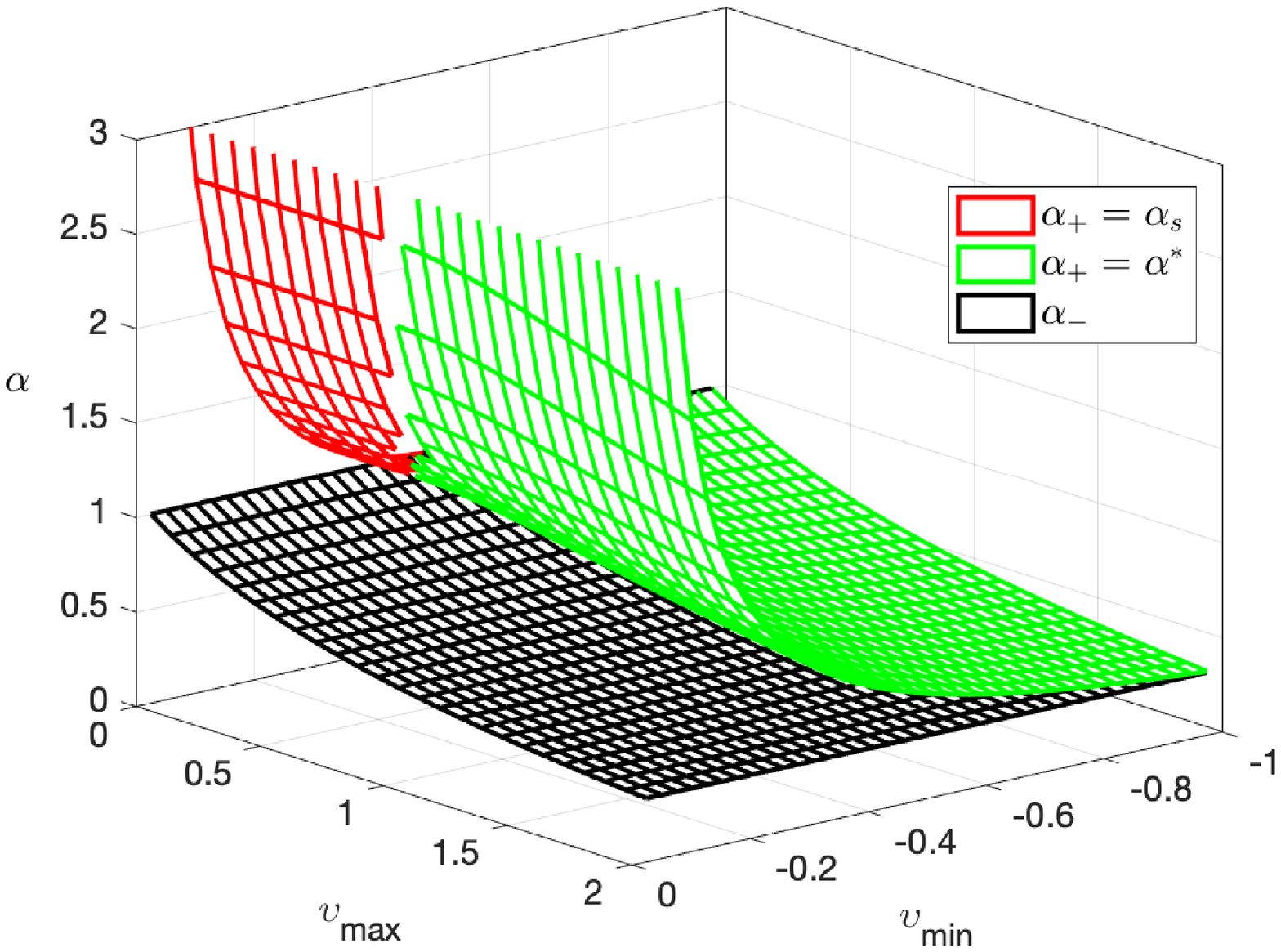}
		\figcaption{Two Critical Thresholds $\alpha_-$ and $\alpha_+$}
		\label{fig: conjectured alpha_max}
	\end{center}

	\vspace{3mm}
	Figure~\ref{fig: conjectured alpha_max} can be used to better understand which triples~$(\alpha,v_{\min},v_{\max}) \in [0,\infty) \times (-1,0) \times (0, \infty)$ lead to all-time positivity. If the triple falls below the surface given by $\alpha_-$, then, according to the Sufficiency Theorem, all-time positivity holds. 
	Alternatively, if the triple lies above the surface given by $\alpha_{+}$, then, according to the Necessity Theorem,  all-time positivity fails. 
	Finally, when the triple lies between the two surfaces, 
we conjecture in Section~\ref{SECTION: Conjs and Supports Of the Conjectures} that all-time positivity also holds.

	
	\vspace{6mm}
	\section{Preliminary Technical Results}\label{SECTION: Technical Results}
	\vspace{0mm} 
	This section provides the technical lemmas underlying the main results.
	In the previous section, we saw that the distinguished path~$v^*$ plays an important role in the motivation of the Necessity Theorem.
	In this section, we provide two preliminary results whose proofs are relegated to Section~\ref{SECTION: Proofs}.
	These preliminaries, involving the behavior of the state $X(v^*,k)$ along the path~$v^*$, are essential to the proof of the main results to follow.
	In addition, later in the paper, these results are seen to provide support for
	our conjecture regarding all-time~positivity. 
	%
	%

	\vspace{5mm}
	\begin{lemma}[Closed-Form for $X(v^*,k)$]\label{Lemma: Closed Form for X(v^*,k)}
		If $\alpha \neq \alpha_s$, for~\mbox{$k\geq 2$}, the state along the distinguished path $v^*$ is given~by 
		\begin{align*}
		X({v^*},k) 
		= \frac{X_0 }{2\sqrt{\theta}}\left( \lambda_+^{k-1}g_{_+}+ \lambda_-^{k-1}g_{_-}  \right) 
		\end{align*} 
		where 
		\begin{align*}
		& \theta  \doteq 4\alpha v_{\min} (1+v_{\min }) + 1,\\[1ex]
		& g_{_{\pm}} \doteq \sqrt{\theta}   \pm \left(2\alpha ({v_{\max }} + 1){v_{\min }} + 1\right)
		\end{align*}
		and
		\begin{align*}
		\lambda_{\pm} \doteq  \frac{1}{2}\left(1\pm{ \sqrt{ \theta } }\right)
		\end{align*}
		are the eigenvalues of $A(v^*,k)$.
		For the singular case,~\mbox{$\alpha = \alpha_s$},
		the state solution is given by
		\begin{align*}
		X({v^*},k) 
		&= \frac{2^{-k} X_0 (( 1-{v_{\max}}+2  {v_{\min}})k+1+{v_{\max}})}{1+v_{\min}}.
		\end{align*} 
	\end{lemma}

	\vspace{5mm}
	\begin{lemma}[Distinguished Path Properties]\label{thm: State Properties along v*} \ \\[1ex]
		(a) If $\alpha > \alpha_s$, then $X(v^*,k)$ is oscillatory about zero and is therefore negative for some values of $k$.
		\\[1ex]
		(b) If $\alpha^* < \alpha \le \alpha_s$ and \mbox{$v_{\max} > 1+2v_{\min}$}, then $X(v^*,k)$ is negative for  sufficiently large $k$. \\[1ex]
		(c)~If $\alpha^* < \alpha \leq \alpha_s$ and~\mbox{$v_{\max} < 1+2v_{\min}$}, then $X(v^*,k)$ is positive for all $k \geq 0$.
		\\[1ex]
		(d) If 
		$0 \leq \alpha \leq \alpha^*$, then~$X(v^*,k)$ is positive for all $k \geq 0$. 
	\end{lemma}

	\vspace{5mm}
	\textbf{Remarks on State Along Distinguished Path $v^*$:} 
	It is interesting in its own right to study the asymptotic behavior of~$X(v^*,k)$, since its tending to zero signifies ``practical bankruptcy,"  even for cases when all-time positivity is assured. 
	For the case \mbox{$0 < \alpha < \alpha_s$} in Lemma~\ref{Lemma: Closed Form for X(v^*,k)}, it is clear that~\mbox{$0<\theta<1$}, which implies 
	$|\lambda_\pm|<1$. Thus, by  the well-known unit-circle stability criterion, for example see~\cite{Jury_1973},~$X(v^*,k) \to 0$ as $k \to \infty$. 
	The closed-form solution for
	the singular case $\alpha = \alpha_s$ also tends to zero, and it is also readily verified, using  l'H\^opital's rule, that the closed-form expression of the state solution~$X(v^*,k)$ is continuous at $\alpha = \alpha_s$. 

	\vspace{6mm}
	\section{Proofs of Preliminary and Main Results}
	\label{SECTION: Proofs}
	\vspace{-1mm}
	This section may be skipped by readers who are not interested in the technical details of the proofs.

	\vspace{5mm}
	\textbf{Proof of Sufficiency Theorem:} Since the case~\mbox{$\alpha =0$} is trivial, we assume $0<\alpha<\alpha_{-}$ and note that 
	it suffices to prove all-time positivity with returns $v(k)$ allowed to range over the larger interval~\mbox{$-1 \leq v(k) \leq v_{\max}$}.  
	We proceed by induction on~$k$. First recall that $\alpha < \alpha_{-}$ was shown to guarantee positivity of $X(2)$ in Section~\ref{SECTION: Main Results}.
	Next, for $k\geq 2,$ we assume~$X(i)>0$ for $i=0,1,\ldots,k$ and all~\mbox{$v(0),\ldots,v(k-1)$}.
	Then for arbitrary \mbox{$v(0),\ldots, v(k)$}, we must show~$X(k+1)>0$.
	Indeed, noting~$1+v(k-1) \geq 0$ and $X(k-1)>0$ by the induction hypothesis, we obtain lower bound
	\begin{align*}
	X(k+1) &= X(k) + \alpha ( 1+ v(k-1) )v(k) X(k-1)\\[.5ex]
	& \geq X(k) - \alpha (1+v(k-1))X(k-1) .
	\end{align*}
	To further lower bound the right hand side above, for~\mbox{$-1 \leq w \leq v_{\max}$}, let $X_{w}(k)$ be the value of $X(k)$ with~$v(k-1)$ replaced by $w$. With this notation, we can write
	$$
	X(k+1) \geq \min_w \bigg\{ X_w(k) - \alpha (1+w)X(k-1) \bigg\}.
	$$
	Since the function to be minimized on the right-hand side above is affine linear in $w$, its minimum value is achieved by~$w=-1$ or~\mbox{$w=v_{\max}$}.
	%
	We now analyze what happens to the minimum in each of case.
	
	\vspace{3mm}
	For $w=-1$, the preceding lower bound of~$X(k+1)$ leads to
	\begin{align*}
	X(k+1) 	
	&\geq X_{-1}(k) >0
	\end{align*}
	by the induction hypothesis. 
	For~\mbox{$w=v_{\max}$}, we~obtain
	\begin{align*}
	X(k+1) 	&\geq  X_{v_{\max}}(k) - \alpha (1+v_{\max})X(k-1).
	\end{align*}
	Since 
	$
	X_{v_{\max}}(k) =  X(k-1) + \alpha ( 1+ v(k-2) ) v_{\max} X(k-2),
	$
	using the facts that $1+v(k-2) \geq 0$, $\alpha >0$, and $X(k-2)$ is positive by the induction hypothesis, it follows that
	\begin{align*}
	X_{v_{\max}}(k) 	&\geq X(k-1) >0
	\end{align*}
	where last inequality holds by induction hypothesis again. Hence, $X(k+1)$ is further lower bounded as
	\begin{align*}
	X(k+1) 	& \geq  [1 - \alpha (1+v_{\max})]X_{v_{\max}}(k).
	\end{align*}
	
	Now, applying the assumed inequality $0<\alpha<\alpha_-$ and the fact that $X_{v_{\max}}(k)>0$ by the induction hypothesis, we obtain~\mbox{$ X(k+1)>0$}.
	\hspace{5mm} $\square$

	\vspace{5mm}
	\textbf{Proof of Lemma~\ref{Lemma: Closed Form for X(v^*,k)}:} Recall the state space representation introduced in Section~\ref{SECTION: INTRODUCTION}. Using the standard state augmentation
	$$
	x(k) \doteq [X(k) \;\; X(k-1)]^T,
	$$
	we obtain the linear time-varying~system
	\[
	x(k+1) = A(v,k)\, x(k)
	\]
	where  $A(v,k)$ is the $2\times 2$ matrix defined in Section~\ref{SECTION: INTRODUCTION}.
	Starting from initial conditions~\mbox{$X(v^*,0)=X(v^*,1) = X_0$} and
	\mbox{$
		X(v^*,2) = (1+\alpha(1+v_{\max})v_{\min})X_0
		$},  in state-space form,
	we have for~\mbox{$k\geq 2$},
	\[
	\begin{bmatrix}
	X(v^*,k+1) \\ X(v^*,k)
	\end{bmatrix} = \underbrace{\begin{bmatrix}
		1 & \alpha (1+v_{\min})v_{\min} \\ 
		1 & 0
		\end{bmatrix}}_{= A(v^*,k)}
	\begin{bmatrix}
	X(v^*,k) \\ X(v^*,k-1)
	\end{bmatrix},
	\]
	and we obtain 
	{\small
		\begin{align*}
		X({v^*},k) = \left[ {\begin{array}{*{20}{c}}
			0&1 
			\end{array}} \right] \begin{bmatrix}
		1 & \alpha (1+v_{\min})v_{\min} \\ 
		1 & 0
		\end{bmatrix}^{k-1}\left[ \begin{gathered}
		X\left(v^*, 2 \right) \hfill \\
		X\left(v^*, 1 \right) \hfill \\ 
		\end{gathered}  \right].
		\end{align*} 
	}We consider two cases: For the generic case, $\alpha \neq \alpha_s$,
	a lengthy but straightforward computation leads to
	\begin{align*} 
	X({v^*},k)  &= \frac{2^{-k}X_0 \left( \left(1+\sqrt{\theta}\right)^{k-1}g_++ \left(1-\sqrt{\theta}\right)^{k-1}g_-  \right)}{\sqrt{\theta}} \\[.5ex]
	&= \frac{X_0 }{2\sqrt{\theta}}\left( \lambda_+^{k-1}g_{_+}+ \lambda_-^{k-1}g_{_-}  \right).
	\end{align*}  
	Another lengthy but straightforward computation shows that~$\lambda_{\pm}$ are  the eigenvalues of $A(v^*,k)$. For the singular case, $\alpha = \alpha_s$,
	we find that
	\begin{align*}
	X({v^*},k) &= \left[ {\begin{array}{*{20}{c}}
		0&1 
		\end{array}} \right]{\left[ {\begin{array}{*{20}{c}}
			1& -1/4 \\ 
			1&0 
			\end{array}} \right]^{k-1}}\left[ \begin{gathered}
	X\left(v^*, 2 \right) \hfill \\
	X\left(v^*, 1 \right) \hfill \\ 
	\end{gathered}  \right]
	\end{align*} 
	which, again, following a third lengthy but straightforward computation, results in
	\[
	X(v^*,k) = \frac{2^{-k} X_0 (k( 1-{v_{\max}}+2  {v_{\min}})+1+{v_{\max}})}{1+v_{\min}}.  \;\;\;\; \square
	\]

	\vspace{5mm}
	\textbf{Proof of Lemma~\ref{thm: State Properties along v*}:} A proof of part (a) that does not use the closed-form of $X(v^*,k)$ can be given immediately by applying Theorem~2.2 in~\cite{Erbe_Zhang_1989}. However, for the sake of self-containment, we provide a first-principles proof here. Assuming that $\alpha > \alpha_s$, we must show the state $X(v^*,k)$ is oscillatory about zero and is negative for some values of $k$.
	By Lemma~\ref{Lemma: Closed Form for X(v^*,k)}, we have
	\begin{align*}
	X({v^*},k) 
	= \frac{X_0 }{2\sqrt{\theta}}\left( \lambda_+^{k-1}g_{_+}+ \lambda_-^{k-1}g_{_-}  \right)
	\end{align*} 
	for~\mbox{$k\geq 2$}.
	With $\alpha > \alpha_s$, it is readily shown that $\theta <0$, which implies that the two eigenvalues $\lambda_{\pm}$ are complex conjugates. 
	It follows that these eigenvalues can be written in polar form as~\mbox{$
		\lambda_{+} = r e^{j \omega}$}
	and~\mbox{$
		\lambda_{-} = r e^{-j \omega}
		$,}
	where $r  = |\lambda_{\pm}|  > 0
	$ and 
	$$
	\omega= \tan^{-1}({\sqrt{|\theta|}}) \in (0,\pi/2).
	$$  
	Next, substituting the polar form of $\lambda_{\pm}$ into $X(v^*,k)$ above,  a lengthy but straightforward calculation shows that
	\[
	X(v^*,k) = B r^{k-1} \cos( (k-1)\omega + \varphi)
	\]
	where $B$ and $\varphi$ are constants, with $B>0$.
	Since $\omega \in (0, \pi/2)$, it is straightforward to find a value of $k$ such that the argument of the cosine lies in $(\pi/2, 3\pi/2),$ thus making the cosine negative. This completes the proof of part~(a).

	\vspace{3mm}
	To prove part (b), we first consider the case $\alpha = \alpha_s$. Then 
	using the formula
	\[
	X(v^*,k) = \frac{2^{-k} X_0 (k( 1-{v_{\max}}+2  {v_{\min}})+1+{v_{\max}})}{1+v_{\min}}
	\]
	for the singular case in Lemma~\ref{Lemma: Closed Form for X(v^*,k)}, for~\mbox{$v_{\max}>1+2v_{\min}$} and $k$ sufficiently large, $X(v^*,k)<0$.  
	Next, for the case~\mbox{$\alpha^* < \alpha < \alpha_s$}, we assume again~\mbox{$v_{\max}>1+2v_{\min}$}.  
	Since~\mbox{$\lambda_+ > \lambda_{-}$}, the state $X(v^*,k)$ will be negative for sufficiently large $k$ if we can show that~\mbox{$g_+ = \sqrt{\theta} + q<0$} where 
	$$
	q \doteq 2\alpha ({v_{\max }} + 1){v_{\min }} + 1.
	$$
	To establish this, since $\alpha \in (\alpha^*, \alpha_s)$, we~have
	$$
	0 < \theta  < 4\alpha^* v_{\min} (1+v_{\min }) + 1 = \frac{(v_{\max}-2 v_{\min}-1)^2}{(1+v_{\max})^2}.
	$$
	Since the square root is an increasing function, the inequality on $\theta$ above implies that 
	$$\sqrt{\theta} < \frac{v_{\max}-2 v_{\min}-1}{1+v_{\max}}.$$
	In addition, we also have
	\[
	q<2\alpha^* (1+v_{\max})v_{\min}+1 = \frac{1-v_{\max}+2v_{\min}}{1+v_{\max}} .
	\]
	Thus, it follows that
	\begin{align*}
	g_+ 
	 	&= \sqrt{\theta} + q \\[.5ex]
	& < \frac{v_{\max}-2 v_{\min}-1}{1+v_{\max}} +  \frac{1-v_{\max}+2v_{\min}}{1+v_{\max}}   =0.
	\end{align*}
	Hence, the proof of part (b) is complete.  
	
	
	\vspace{3mm}
	To prove part (c), we first note that the desired positivity holds trivially for $k=0,1$. For  $k \geq 2$, assuming that $\alpha = \alpha_s$ and $v_{\max} < 1+2v_{\min}$, the singular case formula given in Lemma~\ref{Lemma: Closed Form for X(v^*,k)} leads that
	\begin{align*}
	X(v^*,k) 
	&>  \frac{2^{-k} X_0 (1+{v_{\max}})}{1+v_{\min}}
	\end{align*}
	which is positive for all $k \geq 2$ because $v_{\min}>-1$, $X_0>0$ and $v_{\max}>0$.
	It remains to treat the case $\alpha^* < \alpha < \alpha_s$ and \mbox{$v_{\max} < 1+2v_{\min}$}. 
	To show~\mbox{$X(v^*,k)>0$} for all $k \geq 2,$  
	substitute $g_{\pm} = \sqrt{\theta} \pm q$ and $\lambda_{\pm}= (1\pm \sqrt{\theta})/2$  into~$X(v^*,k)$ and note that $ \theta \in (0,1)$. Then the formula for~$X(v^*,k)$ reduces~to
	{\small \begin{align*}
		X({v^*},k) 
		&= \frac{X_0 }{2^k \sqrt{\theta}} \bigg[ \sqrt{\theta}\bigg(\left(1+{ \sqrt{ \theta } }\right)^{k-1} +\left(1-{ \sqrt{ \theta } }\right)^{k-1} \bigg)\\[.5ex]
		&\hspace{20mm}+   q\bigg(\left(1+{ \sqrt{ \theta } }\right)^{k-1} -\left(1-{ \sqrt{ \theta } }\right)^{k-1} \bigg) \bigg].
		\end{align*}
	}Since
	$v_{\max}<1+2v_{\min}$ and $v_{\min}>-1$, we~obtain	
	$$ 
	q \geq  2\alpha_s ({v_{\max }} + 1){v_{\min }} + 1 =\frac{1-v_{\max}+2v_{\min}}{2(1+v_{\min})}>0. 
	$$
	Since $\sqrt{\theta}>0$, $q >0$ and
	$$
	(1+{ \sqrt{ \theta } })^{k-1} >(1-{ \sqrt{ \theta } })^{k-1}
	$$
	for all $k\geq 2$, it follows that $X(v^*,k)>0.$ This completes the proof of part~(c).
	
	\vspace{3mm}
	Finally, to prove part~(d), since the result trivially follows for~$\alpha = 0$, we assume~$\alpha >0.$ 
	Note that the inequality $\alpha_s \geq \alpha^*$  
	is readily shown to be equivalent to
	$$
	\bigl((1+v_{\max})-2(1+v_{\min})\bigr)^2 \ge 0.
	$$
	Furthermore the above inequalities are both strict
	if and only if $v_{\max} \ne 1+2v_{\min}$. 
	Suppose $v_{\max} \neq 1+2v_{\min}.$ Then~\mbox{$0 < \alpha \leq \alpha^*$} implies~$\alpha < \alpha_s$,
	and so in Lemma~\ref{Lemma: Closed Form for X(v^*,k)}, we have $0<\theta<1$ and~\mbox{$\lambda_{\pm}>0.$} 
	It suffices to prove that~$g_{_\pm}\geq 0$ and that one of $g_+$ or $g_-$ is strictly positive. 
	In the formula for~$g_{_\pm}=\sqrt{\theta} \pm q$,  the quantity~\mbox{$q=
		1+2\alpha(1+v_{\max})v_{\min}
		$} is either negative or nonnegative. 
	If it is nonnegative, then~$g_{_+}>0$, and $g_{_-}\geq 0$ on account of the fact that $\alpha \le \alpha^*$ is equivalent~to
	\[
	\theta \geq [1+2\alpha(1+v_{\max})v_{\min}]^2.
	\]
	Similarly, if the quantity~$q$ above is negative, then $g_{_-} >0,$ while~$g_{+} \geq 0$ on account of the fact that $\alpha \leq \alpha^*$ again.
	
	\vspace{3mm}
	Suppose 
	$
	v_{\max} = 1+2v_{\min}.
	$ 
	Then for the case $\alpha = \alpha^* = \alpha_s$,
	the state~$X(v^*,k)$ for this singular case given in Lemma~\ref{Lemma: Closed Form for X(v^*,k)} applies and is clearly positive for all $k$. 
	Alternatively, for the case $0 < \alpha < \alpha^* = \alpha_s$,
	we argue as in the preceding paragraph and obtain $0<\theta<1$, $\lambda_{\pm}>0$. Moreover, since~$v_{\max}=1+2v_{\min}$, we have $q=\theta$, which leads to
	$$
	g_{_{\pm}} = \sqrt \theta   \pm \theta > 0. 
	$$ This completes the proof of part~(d).
	\hspace{5mm} $\square$
	
	\vspace{5mm}
	\textbf{Proof of Necessity Theorem:}
	Given $\alpha>\alpha_+$, it
	suffices to exhibit a path $v$ for which the state
	$X(v,k)$ is not positive for some~$k$. 
	We claim that the distinguished path $v^*$ is such a path.
	To establish this, we split our analysis into two cases: 
	
	\vspace{3mm}
	\textit{Case 1:} For $v_{\max}\le 1+2 v_{\min}$, we have
	$\alpha_+=\alpha_s$. Thus, it suffices to prove $X(v^*,k)<0$ for some $k$ when $\alpha>\alpha_s$. 
	Using part~(a) of Lemma~\ref{thm: State Properties along v*}, we obtain that the state~$X(v^*,k)$ oscillates and takes negative values for some~$k$.
	
	\vspace{3mm}
	\textit{Case 2:} For \mbox{$v_{\max} >1+2v_{\min}$}, we have $\alpha_+=\alpha^*$. Note that if~\mbox{$\alpha > \alpha_s$}, the negativity of $X(v^*,k)$ is again established by part~(a) of Lemma~\ref{thm: State Properties along v*}. Thus, it suffices to prove~\mbox{$X(v^*,k)<0$} for some $k$ when $\alpha^* < \alpha \le \alpha_s$.
	Since \mbox{$v_{\max} >1+2v_{\min}$}, using part~(b) of Lemma~\ref{thm: State Properties along v*}, we obtain that the state $X(v^*,k)$ is negative for all sufficiently large~$k$. Hence, the proof is~complete. \hspace{5mm} $\square$

	
	\vspace{5mm}
	\section{All-Time Positivity Conjecture and Support}
	\label{SECTION: Conjs and Supports Of the Conjectures}
	\vspace{0mm}
	The conjecture to follow addresses the ``gap" between the lower and upper bounds, $\alpha_{-}$ and $\alpha_{+}$, for all-time positivity provided by the theorems in Section~\ref{SECTION: Main Results}.
	Subsequently, we support the conjecture with analysis and simulations for various cases involving a finite time horizon. As seen below, the notion of ``extreme paths" plays an important~role.

	
	\vspace{5mm}
	{\bf All-Time Positivity Conjecture:}{ \it The all-time positivity condition holds for the gap interval
		$
		\alpha_- \leq \alpha \leq  \alpha_+.
		$}

	\vspace{5mm} 
	{\bf Extreme Paths:} To study the conjecture, for given $N \geq 0$, we consider the~$2^N$ {\it extreme paths}~\mbox{$v^i \in \mathcal{V}^N$}, defined by~\mbox{$v^i(k)$} being either $v_{\min}$ or~$v_{\max}$ for~\mbox{$k=0,1,\ldots, N-1$}. For example, $(v_{\min},v_{\max},v_{\min})$ is an extreme path
	in ${\cal V}^3$.
	First noting that the positivity condition~$X(v,k)>0$ for all $k\leq N$ and all~\mbox{$v\in \mathcal{V}^N$}  is equivalent~to
	$$
	\min_{v \in \mathcal{V}^N } X(v,k) >0
	$$
	for $k \leq N$, 
	we make use of the fact that
	~$X(v,k)$ is
	\textit{multilinear} in~$v$; i.e., affine linear in each component $v(k)$. 
	For example, 
	\[
	X(v,3) = \left[ 1 + v(2)+v(1)v(2) + \alpha\left(v(1)+v(0)v(1)\right) \right]X_0
	\]
	is multilinear in $v(0), v(1)$ and $v(2)$. 
	We now use the well-known fact that the minimum of a multilinear function over a hypercube is attained at one of the vertices; e.g., see~\cite{Barmish_1994}.
	This implies that~$X(v,k)$ is minimized by one of the extreme paths~$v^i$. 
	Hence,~$X(v,k)$ is positive for all~$v \in \mathcal{V}^N$ and all~\mbox{$k \leq N$} if and only~if 
	$$
	\min_{i \in \{1,2,\ldots, 2^N\}} X(v^i,k)>0
	$$ 
	for	~$k=0,1,2,\ldots,N$. 
	For small~$N$, checking this condition is feasible, but for large $N$, the number of ``checks," namely~$2^N$, becomes too large. 
	For example, in the stock market, we can easily have~$N= 100$, but it is computationally prohibitive to check~$2^{100}$ extreme paths.

	\vspace{5mm}
	\textbf{Examples for Various $N$:}
	Taking~\mbox{$X_0=1$}, \mbox{$v_{\max} = 0.9$}, and \mbox{$ v_{\min} = -0.8$}, we have  $v_{\max} > 1+2v_{\min}$, and the gap interval is computed to be~\mbox{$
		[\alpha_-, \, \alpha_+ ] \approx [0.5263, \, 0.5888]. 
		$}
	To support the conjecture, we took $N=10$ and chose~$n=100$ equally-spaced values of~$\alpha$ from the gap interval and, for each~$\alpha$, we used Matlab to check state positivity of each of the~$2^N = 1024$ extreme paths. 
	We found that state positivity held for all of them.  
	In Figure~\ref{fig:State Trajectory_Zoom_In_2},  the~$X(v^i,k)$ are shown for~$
	\alpha = 0.54,
	$
	which lies within the~gap interval above. We also ran many other simulations for various choices of $v_{\min}$, $v_{\max}$ and~\mbox{$N \leq 15$}, and consistently observed that state positivity held in the corresponding gap~interval.

	\vspace{3mm}
	\begin{center}
		\graphicspath{{Figs/}}
		\includegraphics[scale=0.40]{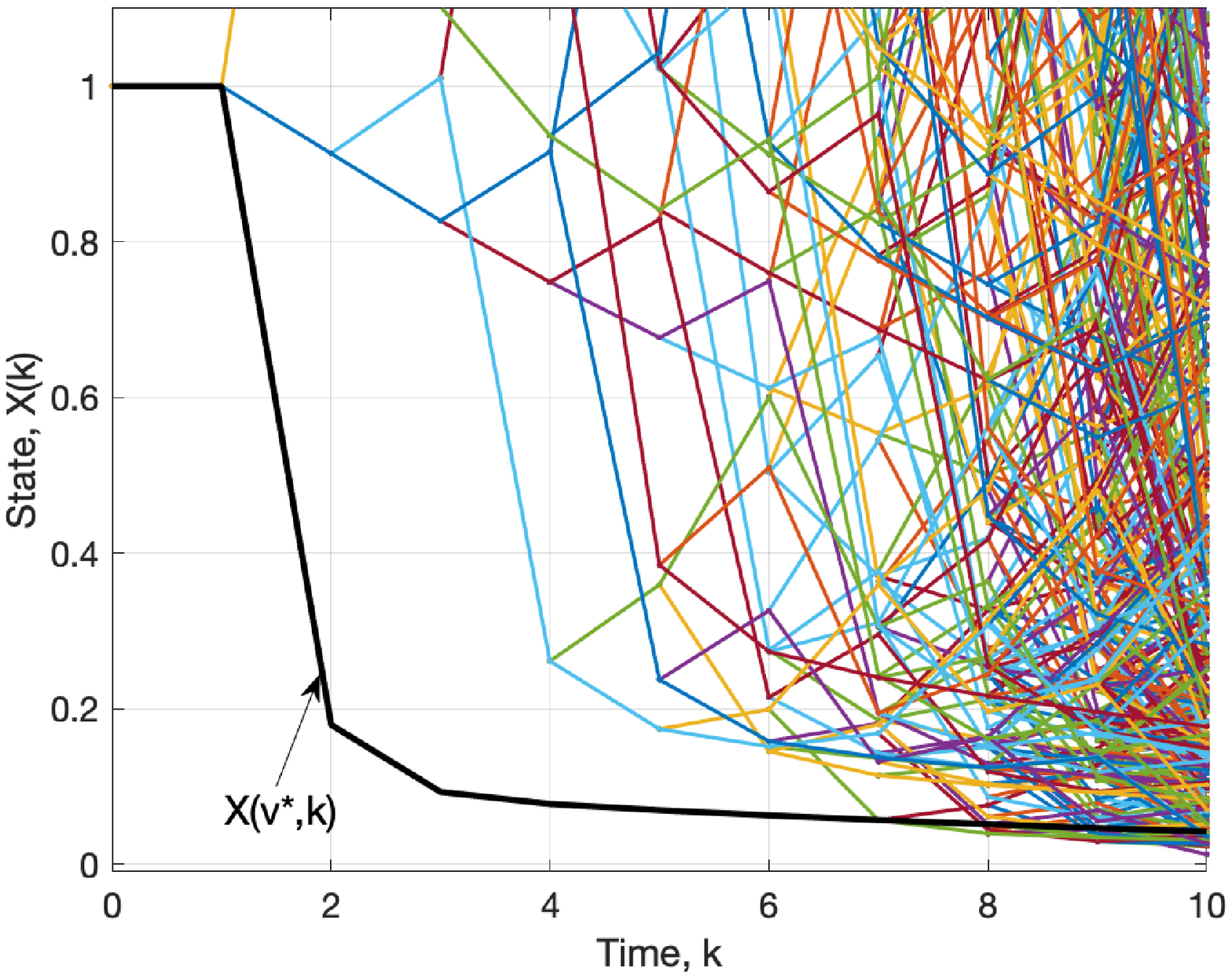}
		\figcaption{Simulation Supporting the Conjecture for $\alpha=0.54$}
		\label{fig:State Trajectory_Zoom_In_2}
	\end{center}
	
	\vspace{3mm}
	Given the motivation for the distinguished path~$v^*$ in terms of  a ``worst-case" trading scenario in Section~\ref{SECTION: Main Results}, it is natural to ask if $X(v^*,k)$ might be the minimum value of $X(v,k)$ 
	for all $k \leq N$.
	However, as seen in Figure~\ref{fig:State Trajectory_Zoom_In_2}, this proves not to be the case for $7\leq k \leq 10$.

	\vspace{3mm}	
	To provide further support for the conjecture, we also studied~\mbox{$N=100$}, \mbox{$X_0 =1$}, \mbox{$v_{\max} = 0.2$} and $v_{\min} = -0.3.$ For~$n=100$ equally spaced values of $\alpha$ in the gap interval~\mbox{$[\alpha_-, \, \alpha_+] \approx [0.8333, \, 1.1905]$}, we generated $200,000$ of the $2^{100}$ extreme paths for each $\alpha$. 
	The positivity condition was seen to be satisfied in all cases.
	Finally, in support of the conjecture, we also ran other simulations for various choices of~$v_{\min}$ and~$v_{\max}$, including smaller values of these bounds to more closely model values found in stock trading, and consistently observed that the desired state positivity held within the gap~interval.

	\vspace{5mm}
	\textbf{Theoretical Result for $N\leq 3$:} In this subsection, we prove that if $\alpha \leq \alpha_+$, then state positivity holds for all partial paths of length $N\leq 3$.
	We begin by noting that the cases $N=0$ and~$N=1$ are immediate since~$X(0)=X(1)=X_0>0$ are the initial conditions.
	Next, for~\mbox{$N=2$}, as shown in the beginning of Section~\ref{SECTION: Main Results}, 
	\begin{align*}
	X(v,2) 
	&\geq [1+\alpha(1+v_{\max})v_{\min}]X_0.
	\end{align*}
	Thus, $X(v,2) >0$ if and only if
	$$
	\alpha < \alpha_{\max}(2) \doteq \frac{1} {|v_{\min}|(1+v_{\max})} .
	$$
	Since it is also easily verified that $\alpha_+ < \alpha_{\max}(2)$, it follows that $X(v,k)>0$ for $v\in \mathcal{V}^2$ and $k \leq 2$ when $\alpha \leq \alpha_+$. 
	The case~\mbox{$N=3$}, per lemma below, requires a lengthier derivation to show that $X(v,3)>0$ if and only if $\alpha < \alpha_{\max}(3)$ where
	$$
	\alpha_{\max}(3) \doteq \frac{1} {|v_{\min}|(2+v_{\max}+v_{\min})} .
	$$
	Then a straightforward calculation shows that \mbox{$\alpha_{+} <\alpha_{\max}(3).$}

	\vspace{5mm}
	\begin{lemma}\label{lemma: Supports of ATP for N=3}
		If  
		$
		\alpha < \alpha_{\max}(3),
		$
		then $X(v,k)>0$ for all~$v\in \mathcal{V}^3$ and all~\mbox{$k \leq 3$}. 
		
	\end{lemma}

	\vspace{3mm}
	\textbf{Proof:} 
	For any~$(v(0),v(1))$, we observe that
	$$
	X(v,3) = X(v,2) + \alpha (1+v(1))v(2) X_0
	$$
	is minimized with $v(2) = v_{\min}$. It follows that
	\begin{align*}
	X(v,3) &\geq X(v,2) + \alpha(1+v(1))v_{\min} X_0\\[1ex]
	& = \left[1 + \alpha \left(\;(1+v_{\min}+v(0))v(1)  +  v_{\min}\right) \right] X_0.	
	\end{align*}
	Since the right-hand side is multilinear in $v(0)$ and $v(1)$, the minimum must occur when they take the values $v_{\min}$ or $v_{\max}$. If $v(1) = v_{\max}$, then to minimize the right-hand side,~$v(0)$ must be $v_{\min}$. 
	In this case, the right-hand side is lower bounded~by
	\[
	\left[1 + \alpha \left(\;(1+v_{\min}+v_{\min})v_{\max} +  v_{\min}\right) \right] X_0.
	\]
	Similarly, if $v(1) = v_{\min},$ then $v(0)$ must be $v_{\max},$ which lower bounds the right-hand side by
	\[
	\left[1 + \alpha \left(\;(1+v_{\min}+v_{\max})v_{\min} +  v_{\min}\right) \right] X_0.
	\]
	It is easy to check that this second bound is strictly smaller than the first. Furthermore, the second bound is positive if and only if $\alpha < \alpha_{\max}(3)$.
	\hspace{5mm}   $\square$

	\vspace{5mm}
	\textbf{Finite-Time Positivity Set:}
	Let $\mathcal{A}(N)$ denote the set of all feedback parameters~$\alpha$ assuring state positivity up to stage~$N$. Define 
	$$
	\alpha_{\max}(N) \doteq \sup \{\alpha \geq 0: [0,\alpha) \subseteq \mathcal{A}(N) \}.
	$$ 
	Then we have already seen above that
	$
	{\cal A}(2) = \left[ 0 ,\, \alpha_{\max}(2) \right)
	$
	and~\mbox{$
		{\cal A}(3) = \left[ 0 , \,
		\alpha_{\max}(3) \right)
		$}
	with $\alpha_{\max}(3) < \alpha_{\max}(2)$ readily verified. Beyond these two simple cases,
	one can in principle determine whether or not a given feedback
	parameter $\alpha$
	belongs to~${\cal A}(N)$ by checking all extreme paths.
	We also know, by
	the Sufficiency Theorem, that~\mbox{$
		[0, \, \alpha_-) \subseteq \mathcal{A}(N)
		$}. If the All-Time Positivity Conjecture is true, we must have~\mbox{$[0, \, \alpha_+] \subseteq \mathcal{A}(N)$} as well. 
	Moreover, since~\mbox{${\cal A}(N+1) \subseteq {\cal A}(N)$} for all $N$,  the $\alpha_{\max}(N)$ are nonincreasing, and since they are bounded below by $\alpha_-$, they converge to a limit 
	$$
	\alpha_{\infty} \doteq \lim_{N\to\infty} \alpha_{\max}(N). 
	$$ 
	It is also readily verified that
	$
	\alpha_{\infty} \leq \alpha_+;  
	$ otherwise, there would exist an $\alpha \in (\alpha_+, \alpha_\infty)$ assuring all-time positivity, which contradicts the Necessity Theorem.
	Finally, if the All-Time Positivity Conjecture is true, then 
	$\alpha_{\infty} \geq \alpha_+,
	$
	in which case it would follow that~\mbox{$\alpha_{\infty} = \alpha_+.
		$ }
	%
	%
	
	\vspace{6mm}
	\section{Conclusion and Future Work}
	\label{SECTION: CONCLUSION AND FUTURE WORK}
	\vspace{0mm}
	In this paper, we considered a state positivity problem motivated by trading risky assets in the presence of delay.
	The desired positivity of the state was studied in terms of two critical thresholds,~$\alpha_-$ and  $\alpha_+$ with~\mbox{$\alpha_- < \alpha_+$}.
	First we proved that~\mbox{$\alpha < \alpha_-$} is sufficient for all-time positivity. 
	Then we  proved that~\mbox{$\alpha > \alpha_+$} is necessary for all-time positivity. 
	Finally, we conjectured that state positivity is guaranteed for the ``gap" interval~\mbox{$\alpha_- \leq \alpha \leq \alpha_+.$}
	Support for this conjecture, both theoretical and computational, was also~provided.
	


	\vspace{3mm}
	Regarding further research, we mention two attractive directions. 
	The first is obviously to pursue a proof  of the conjecture.    
	Based on many simulations, we consistently observed the following phenomenon: \textit{If~\mbox{$X(v^*,k)>0$} for~\mbox{$k\leq N$}, it follows that $X(v,k)>0$ for~$k\leq N$ and all~\mbox{$v\in \mathcal{V}^N$}}; e.g., see Figure~\ref{fig:State Trajectory_Zoom_In_2} where $X(v^*,k)$ is positive and the other states $X(v^i,k)$ are positive too, which implies~$X(v,k)>0$ for $k\leq 10$. 
	If this observation is true for all $N$, then parts (c) and (d) of Lemma~\ref{thm: State Properties along v*} give us all-time positivity for~\mbox{$\alpha \leq \alpha_+$}.

	\vspace{3mm}
	A second direction for future research involves studying the state positivity problem when~$v(k)$ is vector-valued rather than a scalar. That is, if~\mbox{$v(k) \in \mathbb{R}^m$} with~$v_i(k)$ being the~\mbox{$i$th} component satisfying~\mbox{$v_{\min,i} \leq v_i(k) \leq v_{\max,i}$} with \mbox{$-1<v_{\min,i}<0<v_{\max,i}$} for $i=1,2,\ldots,m$, then, motivated by portfolio rebalancing problems with delay,  the more general state equation
	\[
	X(k+1) = X(k) + \sum_{i=1}^{m} \alpha_i (1+v_i(k-1))v_i(k)X(k-1)
	\]
	arises where the $\alpha_i \geq 0$ are scalar constant feedback parameters. 
	In this case,  generalization of the theory in this paper would be of~interest. 
	To this end, one result along these lines is that the condition
	$$
	4\sum_{i=1}^m \alpha_i (1+v_{{\min},i})|v_{\min,i}|  > 1,
	$$
	leads to oscillation and failure of all-time positivity. This can be established using arguments similar to those given in the proof of Lemma~\ref{thm: State Properties along v*} and the related literature.




\end{document}